\documentclass[12pt]{article} 
\usepackage{amsfonts,amssymb} 
\usepackage[cp1251]{inputenc}

\usepackage[margin=25mm, a4paper]{geometry}  

\usepackage{tikz}

\def\n{\noindent}  \def\?#1{}
\def\IZ{{\mathbb{Z}}}  \def\IR{{\mathbb{R}}}
 \def\cA{{\cal A}}   \def\cM{{\cal M}}

   \def\phi{\varphi}  
\def\ep{\varepsilon}  \def\t{\tilde} \def\gD{{\mathbb{D}}}
     \def\v#1{\overline{#1}} 

\def\mod1{\,({\rm mod\ } 1)\,} 
  \def\bdelete#1{} 

\def\beq#1#2{\begin{equation} \label{#1} #2 \end{equation}}
\def\bea#1{\begin{eqnarray*} #1 \end{eqnarray*}} \def\a{\!\!\!&\!\!\!\!&}

\def\toas#1{\stackrel{#1}{\longrightarrow}}

\def\proof{\smallskip \noindent {\bf Proof. \ }}       
\def\blanksquare{\,\,\,$\sqcup\!\!\!\!\sqcap$}         
\def\qed{\hfill\blanksquare\linebreak\smallskip\par}   
\def\hence{\Longrightarrow}                                 

\def\thname{Theorem}  \def\lmname{Lemma}    \def\prname{Proposition}
\def\dfname{Definition}  \def\crname{Corollary}  \def\rmname{Remark}
\def\exname{Example}   \def\conname{Hypothesis}

\newtheorem{theorem}{\thname}[section]   
\newtheorem{lemma}{\lmname}[section]     
\newtheorem{proposition}[lemma]{\prname} 

\newtheorem{dftn}{\dfname}[section]
\newenvironment{definition}{\begin{dftn}\rm}{\end{dftn}} 
\def\bdef#1{\begin{definition} #1 \end{definition}}
\newtheorem{rmrk}[lemma]{\rmname}
\newenvironment{remark}{\begin{rmrk}\rm}{\end{rmrk}}     

\makeatletter\def\fps@figure{htbp}\makeatother 

\begin{document}

\title{Local against global dynamical entropies}
\author{Michael Blank\thanks{
        Moscow Center for Advanced Studies;}
        \thanks{National Research University ``Higher School of Economics''; e-mail: mlblank@gmail.com} 
          ~and Maxim Polyakov\thanks{National Research University ``Higher School of Economics''; e-mail: pmaxol73@gmail.com}}
\date{\today} 
\maketitle

\begin{abstract} 
We study the relationship between the ``complexity'' of a dynamical system, defined 
by the classical Kolmogorov-Sinai metric entropy construction, and the ``complexity'' 
of its typical trajectories, defined by M. Blank's local dynamic entropy. 
A key feature is the narrowing (compared to the Kolmogorov-Sinai construction) 
of the choice of the family of admissible space partitions, which turns out to be related 
to the statistics of a specific trajectory. This approach is extended to the analysis 
of arbitrary sequences (which may not be included in any dynamical system).
\end{abstract}

{\small\n
2020 Mathematics Subject Classification. Primary: 37A44; Secondary: 37A35, 11N05, 11K65.\\
Keywords: ergodic theorem, dynamical system, metric entropy, topological entropy, complexity.
}

\tableofcontents%

\section{Introduction}\label{s:intro} 

The now-classical constructions of the Kolmogorov-Sinai metric entropy and 
the Adler-Cohnheim-McAndrew topological entropy provide a fairly complete 
quantitative characterization of the ``complexity'' of a dynamical system as a whole.
However, these approaches provide no information about the "complexity" of individual
trajectories of this system, let alone the "complexity" of sequences
that are not trajectories of any system. Recently, M. Blank (2025) proposed
an original construction of local dynamical entropy to solve this problem
and applied it to the analysis of the ``complexity'' of various number-theoretic sequences.
As in the Kolmogorov-Sinai construction, local entropy is determined in two steps: first, 
the conditional entropy is determined for a given measurable partition of space, and 
then the supremum over the set of admissible partitions is taken to determine the unconditional entropy.

In the case of number-theoretic sequences, we are dealing with sequences
with a finite alphabet, and the second step is trivial. In the general case, however,
we will show that the choice of the set of admissible partitions (unlike in the Kolmogorov-Sinai 
and Adler-Cohnheim-McAndrew constructions) has a fundamental impact on the result. 
In particular, we will show that when choosing admissible partitions in the Kolmogorov-Sinai 
style, sequences that would naturally be considered ``simple'' (for example, $\{1/n\}$) 
will have infinite entropy (which contradicts common sense).

In this paper, we will clarify the relationship between the Kolmogorov-Sinai metric entropy
and the dynamical local entropy of its typical trajectories. This will first be done for conditional 
entropies over arbitrary measurable partitions, and then for unconditional ones. It is in this 
second stage that the correct choice of the set of admissible partitions will be determined, 
which will allow us to extend this construction to arbitrary sequences (which may not be 
included in any dynamical system). 
To distinguish between the dynamical entropies describing the ``complexity''  of 
a dynamical system as a whole and entropies corresponding to the ``complexity'' of 
individual sequences/trajectories, we use the word "global" in the former case and "local" in the latter. 
The main result is that the local entropy of a typical trajectory of an ergodic system is 
equal to the metric entropy of the system (Theorem~\ref{t:global-erg}), 
while in the non-ergodic case the metric entropy is equal to the integral over local 
entropies of all trajectories  (Theorem~\ref{t:global-nonerg}).

This work required solving a number of complex technical problems,
from developing a new approach to studying the properties of unconditional local
entropy to analyzing classical metric entropy, where it was necessary
to prove that the supremum over all measurable partitions coincides with the supremum over
the aforementioned admissible partitions, which depend on specific trajectories.

The paper is organized as follows. In Section~\ref{s:Entropy}, we recall
the classical definitions associated with the concept of Shannon entropy and 
classical dynamical entropies, as well as local entropies of individual sequences. 
Section~\ref{s:glob-loc-cond} is dedicated to the comparison between 
metric and local conditional entropies, while in Section~\ref{s:glob-loc} 
we study the same question for unconditional entropies. In the latter case we need 
to make additional assumptions about continuity of dynamical systems under study. 
A large part of Section~\ref{s:glob-loc} is devoted to the choice of the set 
of admissible partitions and their properties. In Section~\ref{s:top-Inf} we discuss 
connections between topological and information entropies. 
Finally, in Section~\ref{s:topology} we prove the topological invariance of local dynamical entropies.

\section{Entropy type constructions in ergodic theory}\label{s:Entropy}
Let $(X,\Sigma)$ be a measurable space with a Borel $\sigma$-algebra $\Sigma$ of measurable subsets. 

\bdef{Iterations of a measurable map $f:X\to X$ define a discrete time dynamical system -- a pair $(X,f)$.
A trajectory of this system starting at a point $x_0\in X$ is a sequence 
$\v{x}:=\{x_n:=f^nx_0\}_{n=0}^\infty$.}

Here $f^nx:=\underbrace{f\circ f\circ\dots\circ f}_{n} x$ is the $n$-th iteration of the map $f$ at a point $x$.

\bdef{A $\sigma$-additive measure $\mu:\Sigma\to\IR_+$ is said to be $f$-invariant if 
$\mu(f^{-1}A) = \mu(A) ~~\forall A\in\Sigma$.}

Recall that an $f$-invariant measure $\mu$ is called {\em ergodic} if 
for each $f$-invariant measurable set $A\in\Sigma$ either $\mu(A)=0$ 
or $\mu(X\setminus A)=0$.

In this section we recall basic definitions and properties of entropy type constructions 
related to dynamical systems. We start with the classical objects: static (Shannon) entropy 
of a discrete distribution and three dynamical global ones - Kolmogorov-Sinai metric, 
Adler-Konheim-McAndrew topological and Brin-Katok local entropies. Then we move to 
local and information entropies of individual sequences/trajectories introduced in a recent 
paper by the first author \cite{Bl25}.

\subsection{Entropy of a discrete measure}\label{s:sEnt}
Let $\Omega:=\{\omega_1,\omega_2,\dots\}$ be a discrete set, equipped with a finite 
(necessarily probabilistic) distribution $\v{p}:=(p_1,p_2,\dots)$ (i.e. $p_i\ge0,~\sum_ip_i<\infty$).
A distribution $\v{p}$ is said to be probabilistic if $||\v{p}||:=\sum_ip_i=\v{p}(\Omega)$.

\bdef{For a nonnegative number $t$ we define a function $H(t):=-t\log t$ 
(to which we will refer as an entropy), where $\log(\cdot):=\log_2(\cdot)$, 
and by continuity we set $H(0)=0$.}

\bdef{The (Shannon) {\em entropy} of a distribution $\v{p}$ is defined as $H(\v{p}):=\sum_i H(p_i)$.}

To distinguish this notion from the ones considered for dynamical systems we will refer to 
this entropy as static.

\subsection{Classical dynamical entropies}\label{s:dEnt}
Now we are ready to introduce and discuss in some detail classical approaches to the notion 
of entropy in ergodic theory which we will need for comparison with with our own notions 
of local entropies.

We start with the classical Kolmogorov-Sinai construction (see, e.g.\cite{Wa, Bi}, for detail). 

Recall that a finite collection of sets $\Delta:=\{\Delta_i\}$ is called a measurable partition\footnote{
        In what follows, if the measure $\mu$ is not chosen, the assumption 
        $\mu(X \setminus (\cup_i \Delta_i ))=0$ is dropped.} 
of $(X,\Sigma,\mu)$ if $\mu(X \setminus (\cup_i \Delta_i ))=0$ and $\Delta_i \cap \Delta_j=\emptyset~~\forall i\ne j$.

\bdef{A {\em refinement} of a pair of measurable partitions $\Delta, \Delta'$ of $(X,\Sigma)$ 
is a new partition $\Delta \bigvee \Delta' := \{\Delta_{i}\cap \Delta'_{j}: ~~\forall i,j\}$.}

Fix a dynamical system $(X,f,\Sigma,\mu)$ with a finite $f$-invariant measure $\mu$. 

\bdef{By refining the union of all pre-images under the action of the map $f$ of the partition $\Delta$,  
i.e. $\{f^{-1}\Delta_{i}\}$, we get a new measurable partition, which we denote by $f^{-1}\Delta$. 
Under the  {\em $n$-th refinement} of the partition $\Delta$ one means
$$\Delta^n:=\Delta^{n-1}\bigvee f^{-1}\Delta^{n-1}, ~~\Delta^0:=\Delta .$$}

\bdef{By a {\em conditional Kolmogorov-Sinai entropy} of the map $f$ with respect to the 
partition $\Delta$ one means 
$$h_\mu(f|\Delta):=\lim\limits_{n\to\infty}\frac1n H_\mu(\Delta^n),$$ 
where $H_\mu(\Delta):=-\sum_{i}\mu(\Delta_{i})\log\mu(\Delta_{i})$ is the entropy of 
the discrete distribution $\{\mu(\Delta_{i})\}$.}

Observe that after the $n$-th refinement we get a finite partition and the collection of 
numbers $\{\mu(\Delta_j^n)\}_j$ is a finite distribution. Therefore $H_\mu(\Delta^n)$ 
is the (static) entropy of this distribution. Note also that when defining the partition 
refinement, elements of measure zero are often excluded. This is not necessary since 
$H(0)=0$ and thus elements of measure zero do not contribute to the dynamical entropy.

\bdef{By a {\em metric Kolmogorov-Sinai entropy} of the dynamical system $(X,f,\Sigma,\mu)$ 
one means $h_\mu(f):=\sup_\Delta h_\mu(f|\Delta)$, 
where the supremum is taken over all finite measurable partitions $\Delta$.}

If the map $f$ is continuous R.L. Adler, A.G. Konheim, M.H. McAndrew proposed 
in \cite{AKM} another approach for the construction of entropy type quantities, 
based on covers rather than partitions without connections with invariant measures.

Let $\Delta:=\{\Delta_i\}_1^r$ be a finite {\em covering} of $X$ by open sets.
Denote by $A_\Delta^n$ the smallest subcover of the $n$-th refinement of $\Delta$. 
Note that in general the collection of open sets $\Delta^n$ needs not to be a cover of $X$.

Using this notation one defines the conditional topological entropy as follows  
$$h_{{\rm top}}(f|\Delta):=\lim_{n\to\infty}\frac1n \log(\#A_\Delta^n) .$$

\bdef{The {\em topological entropy} is defined as 
$h_{{\rm top}}(f):=\sup_\Delta h_{{\rm top}}(f|\Delta)$.}

A few comments: 
$h_\mu(f)$ sensitively depends on the $f$-invariant measure $\mu$, 
while $h_{{\rm top}}(f)$  makes sense only for $C^0$-maps. 
Additionally, a measurable dynamical system need not admit invariant measures 
(not speaking about its continuity).

There is a number of other approaches to define global entropy of a system, among 
which we highlight a version proposed by Brin and Katok \cite{BK}, allowing to 
study the complexity of an individual trajectory of a continuous map. 

Let $f$ be a continuous map from  a compact metric space $(X,\rho)$ into itself. 
To take into account the correspondence between trajectories starting at different 
points one introduces yet another metric.

\bdef{The $n$-th {\em Bowen metric} $\rho_n$ on $X$ between the points $u$ and $v$ is    %
$$\rho_n(u, v):=\max\left\{\rho\left(f^k(u), f^k(v)\right):\;k =0, \ldots,n-1 \right\}.$$}
This quantity controls the divergence between first $n$ points of trajectories starting at $u$ and $v$. 

Denote by $B_\ep^n(u)$ the open ball of radius $\ep$ in the metric $\rho_n$ around $u\in X$.

\bdef{{\em Brin-Katok local entropy} \cite{BK} of $f$-invariant measure $\mu$ 
at point $u\in X$ is 
$h_{BK}(f,\mu,u):=-\lim\limits_{\ep\to0}
              \lim\limits_{n\to\infty}\frac{1}{n}\log\mu(B_\ep^n(u))$.}%

An important point here is that this quantity depends heavily upon the behavior of all 
neighboring trajectories, which excludes the generalizations of interest to us for the 
analysis of a separate sequence not included in a specific dynamical system.

\subsection{Simple properties of global entropies} 
Let $\sigma(\Delta)$ be the minimal $\sigma$-algebra generated by the sets $\{f^{-n}\Delta_i\}_{n\in\IZ_+}$.

\bdef{A partition $\Delta$ of $X$ is called {\em generating} if $\sigma(\Delta)=\Sigma$.  }

\begin{proposition}\label{p:entropy} \cite{Bi,EW,Wa,BK}
\begin{enumerate}
\item If $\sigma(\Delta)=\Sigma$, then $h_\mu(f|\Delta)=h_\mu(f)$
\item Let $\Delta^{n}$ be an increasing sequence of finite partitions, such that 
        $\sigma(\Delta^{n})\toas{n\to\infty}\Sigma$. Then $h_\mu(f|\Delta^{n})\toas{n\to\infty}h_\mu(f)$
\item $h_\mu(f^n)=nh_\mu(f)~\forall n\in\IZ_+$. 
\item $h_\mu(f)$ is preserved under the action of an isomorphism.
\item  $h_{{\rm top}}(f) = \sup_\mu h_\mu(f)$ if $f$ is continuous and supremum is taken 
         over all $f$-invariant measures.

\item if $\mu$ is ergodic, then $h_\mu(f) = h_{BK}(f,\mu,u)$ for $\mu$--a.a. $u\in X$, 
         otherwise $h_\mu(f) = \int h_{BK}(f,\mu,u) ~d\mu(u)$. 
\end{enumerate}
\end{proposition}

The restriction to continuous maps is rather severe. Therefore let us try to emulate the 
construction of the topological entropy for a measurable map by a partition $\Delta$ instead of a covering. 
It turns out that in this generalized setting one gets only 
$$ \hat{h}_{{\rm top}}(f) \ge \sup_\mu h_\mu(f) $$ 
and the this inequality may be strict.

\subsection{Local dynamical entropies}\label{s:locdyn}
Global entropies characterize the complexity of the dynamic system as a whole, but do not provide any information about the complexity of an individual trajectory (except for the Brin-Katok version, 
which heavily depends on the behavior of trajectories nearby).  

We therefore discuss two completely different constructions of local entropies (introduced in \cite{Bl25})  
that are based only on the statistical properties of an individual sequence, and study their connections 
with global entropies in the case where this sequence is indeed a true trajectory.

Let $\Delta:=\{\Delta_i\}_{i=1}^r$ be a finite measurable partition of $(X,\Sigma)$. 
We identify indices of the elements of the partition $\Delta$ with letters of an alphabet 
$\cA:=\{a_i\}_{i=1}^r$, which will allow us to further encode arbitrary sequences of points 
from $X$ with these letters as follows. We encode a given sequence $\v{x}$ with a sequence 
of letters from the alphabet $\cA$, corresponding to the partition elements that
the sequence successively visits, and study the statistics of subwords in the initial
segments of this encoding.

$$ \v{x} = (x_1x_2\dots x_k\underbrace{x_{k+1}\dots x_{k+n}}_{\v{w}} 
                 x_{k+n+1}\dots  x_N\dots) $$

\bdef{For a given sequence $\v{x}:=(x_1,x_2,\dots)$ we say that its initial 
segment of length $N$ contains a word $\v{w}:=(w_1,\dots,w_n)$, formed by the 
letters $w_i\in\cA$, if there is an index $i$ such that
$$x_{i+j}\in \Delta_{w_j} ~~\forall j<n+1, i+j \le N .$$}

Denote by $\v{p}(\v{x},n,N)=(p_1,p_2,\dots)$ the distribution (set of frequencies) of 
various words of length $n$ in the starting segment of length $N$ of the 
sequence $\v{x}:=(x_1,x_2,\dots)$, and by $L(\v{x},n,N)$ the number of 
different words $\v{w}$ of length $n$ in this segment.

Applying Shannon entropy to characterize the distribution $\v{p}(\v{x},n,N)$ and passing 
to the limit we come to the following notion.

\bdef{{\em The conditional local entropy} of the sequence $\v{x}$ with respect to the partition
$\Delta$ is %
\beq{e:B-ent-loc-cond}{ 
   h_{{\rm loc}}(\v{x}|\Delta) 
            := \limsup_{n\to\infty} \limsup_{N\to\infty} \frac1n H(\v{p}(\v{x},n,N)) ,}}    
where 
$$ H(\v{p}(\v{x},n,N)) := -\sum_{i=1} p_i\log p_i $$
is the entropy of the finite discrete probability distribution $\v{p}(\v{x},n,N)$.

\bdef{{\em The local entropy} of the sequence $\v{x}$ is 
\beq{e:B-ent-loc}{ h_{{\rm loc}}(\v{x}) := \sup_{\Delta\in\gD} h_{{\rm loc}}(\v{x}|\Delta),}   
where the supremum is taken over admissible partitions $\Delta\in\gD$.}

\bdef{By the family of {\em admissible} partitions $\gD$ we mean all finite
measurable partitions $\Delta:=\{\Delta_i\}_{i=1}^r$ of the measurable space $(X,\Sigma)$,
the elements of which satisfy some additional conditions, which will
be introduced depending on the context.}

In fact, defining Kolmogorov-Sinai metric entropy, not all finite partitions are taken into account, but only 
partitions with elements of positive $\mu$-measure. In section~\ref{s:glob-loc} we discuss how 
various admissibility assumptions influence the local entropies.

\begin{remark}
In a sense, our local construction follows the same logic as the global construction 
of the Kolmogorov-Sinai metric entropy, but uses only information about the given sequence.
\end{remark}

\bigskip

Similarly to the above construction, one implements a local analogue of the topological entropy.

For a given admissible partition $\Delta$, denote by $L(\v{x},n,N)$ the number of  
{\em different} words of length $n$ in the starting 
segment of length $N$ of the sequence $\v{x}:=(x_1,x_2,\dots)$

\bdef{{\em The conditional information entropy} of the sequence $\v{x}$ with respect to 
the partition $\Delta$ is %
\beq{e:B-ent-inf-cond}{ h_{{\rm info}}(\v{x}|\Delta) 
   := \limsup_{n\to\infty} \limsup_{N\to\infty} \frac1n\log L(\v{x},n,N) .}}   

\bdef{{\em The information entropy} of the sequence $\v{x}$ is 
\beq{e:B-ent-inf}{ h_{{\rm info}}(\v{x}) := \sup_{\Delta\in\gD} h_{{\rm info}}(\v{x}|\Delta) ,}    
where the supremum is taken over admissible partitions $\Delta\in\gD$.}

\subsection{Connections between local and information entropies}\label{s:Loc-Inf}
It is known (see, for example, \cite{Ka}) that under sufficiently weak constraints imposed 
on the dynamical system $(f,X,\Sigma)$, the relation
$$ h_{{\rm top}}(f) = \sup_\mu h_\mu(f) $$ is satisfied. Here one maximizes metric 
entropies over all invariant measures. Despite that in our local entropy construction 
we have no such free parameter, a similar but weaker relation can be established 
between local and information entropies.

\begin{lemma}\label{l:ineq-entr} 
For each sequence $\v{x}$ and each finite measurable partition $\Delta$ we have  
$$ h_{{\rm loc}}(\v{x}|\Delta) \le h_{{\rm info}}(\v{x}|\Delta) .$$   
\end{lemma}

\proof For a given finite partition $\Delta$, an arbitrary sequence
$\v{x}$ and arbitrary $n<N$, consider the discrete distribution
$\v{p}(\v{x},n,N)=(p_1,p_2,\dots)$ introduced in the previous section. 

Since $L(\v{x},n,N)$ is the number of distinct words of length $n$ in the starting segment 
of length $N$ of the sequence $\v{x}$, we observe that at most $L(\v{x},n,N)$
elements of the discrete distribution $\v{p}(\v{x},n,N)$ are positive.
Now using the fact that the Shannon entropy of a finite probability
distribution reaches its maximum on the uniform distribution, we obtain
$$ H(\v{p}(\v{x},n,N)) := -\sum_{i=1} p_i\log p_i \le \log L(\v{x},n,N) ,$$
whence the desired assertion follows. \qed

Questions related to unconditional entropies will be discussed later in Section~\ref{s:glob-loc}, 
after the introduction of the proper family of admissible partitions. 

\bigskip


\begin{lemma}\label{l:per} If the sequence $\v{x}$ is (eventually) periodic,
then its distribution -- the measure $\mu_{\v{x}}$ is well defined,
and $h_{{\rm loc}}(\v{x}) = h_{{\rm info}}(\v{x}) = 0$.
\end{lemma}
\proof We start with the periodic case with a certain period $\ell<\infty$, 
where the existence of $\mu_{\v{x}}$ is evident.
For each $n>\ell$ in $\v{x}$ there are at most $\ell$ distinct words of length $n$,
each of which has frequency $1/\ell+o(1/N)$. 
By definition, the local entropy of this sequence is $\lim_{n\to\infty}\frac1n\log\ell=0$.
Similarly, $L(\v{x},n,N)\le\ell~~\forall (n,N)$, from which we obtain the zero value of
information entropy.

In the eventually periodic case the sequence becomes periodic only starting from 
some finite index $j$. On the other hand, the limit relations 
(\ref{e:B-ent-loc-cond}, \ref{e:B-ent-loc}, \ref{e:B-ent-inf-cond}, \ref{e:B-ent-inf}) 
do not change if we discard a finite starting segment of a sequence. Therefore 
one applies the argument above for the sequence starting from the index $j$ 
to get the result. \qed 


\section{Connections between metric and local conditional entropies}\label{s:glob-loc-cond}
\subsection{Ergodic case} \label{s:glob-loc-cond-erg}

Since the sequences $\v{x}:=\{x_k\}$ in this Section will be trajectories of a certain 
dynamical system defined by the action of a map $f:X\to X$, they are defined completely 
by initial points. Therefore to simplify notation for the trajectory  $\v{x}:=\{x_k:=f^kx\}$ 
we denote by $x$ this initial point. 

\begin{theorem}\label{t:glob-loc-erg} 
Let $\mu$ be an ergodic $f$-invariant measure. Then 
$h_\mu(f|\Delta) = h_{{\rm loc}}(\v{x}|\Delta)$ 
for each finite measurable partition $\Delta$ and for $\mu$-a.a. $x\in X$. 
\end{theorem}
\proof 
Observe that for a given map $f$ and a partition $\Delta$ its $n$-th refinement 
$\Delta^n:=\{\Delta_i^n\}$ corresponds to $n$-letter words in the alphabet generated 
by the elements of $\Delta$. 

Before we give the proof, let's discuss a few pitfalls associated with 
the ergodic constructions we use.
To estimate the frequency of various words in the ``coding'' of $\v{x}$, for example, the 
two-letter word $ij$, we must calculate the following sum:
$$ \frac1n\sum_{k=0}^{n-1}1_{\Delta_i}(f^kx) 1_{\Delta_j}(f^{k+1}x) .$$
This can be written in terms of ergodic averages 
$$ \frac1n\sum_{k=0}^{n-1}1_{\Delta_i\times\Delta_j}(F^k(x,y)) $$
along the trajectory of the mapping $F:=f\times f:X\times X \to X\times X$,
which starts at the point $(x,fx)$. Assuming ergodicity of the direct 
product system $(F, X\times X, \mu\times\mu)$, we obtain the result
using Birkhoff's classical ergodic theorem. The only problem here is that
despite the ergodicity of the system $(f, X, \mu)$, the direct product system
needs not be ergodic. So instead we will follow a different idea, based on 
the observation that the word $ij$ appears if and only if
$f^kx\in \Delta_i \cap f^{-1}\Delta_j$ for some $k$. 
Consequently, we consider another ergodic mean
$$ \frac1n\sum_{k=0}^{n-1}1_{\Delta_i \cap f^{-1}\Delta_j}(f^kx) ,$$
for convergence of which we only need the ergodicity of the original system.

Now we are ready to give the formal proof. 
By the definition of the conditional Kolmogorov-Sinai entropy
$$ h_\mu(f|\Delta) = \lim_{n\to\infty}\frac1n H_\mu(\Delta^n) 
                            = \lim_{n\to\infty}\frac1n \sum_{\Delta_i^n\in\Delta^n} H(\mu(\Delta_i^n)). $$ 
On the other hand, by Birkhoff ergodic theorem (see e.g. \cite{Bi}) $\forall A\in\Sigma$ 
and $\mu$-a.a. $z\in X$ %
\beq{e:Birkhoff}{ \lim_{N\to\infty} \frac1N \sum_{k=0}^{N-1} 1_A(f^kz) = \mu(A).}

Let us show that the passage to the limit and the summation on the left side of this formula 
can be interchanged. The relation (\ref{e:Birkhoff}) implies that $\forall \delta>0~\exists N'$  
$$ \left|\frac1N \sum_{k=0}^{N-1} 1_A(f^kz) - \mu(A)\right| < \delta~~ \forall N>N' .$$
Using continuity of the function $H(\cdot)$ in the segment $[0,1]$ we that 
$$ \forall z\in[0,1]~~\forall \ep>0~\exists \delta(z)>0:~~ |z-y|<\delta(z) \hence |H(z) - H(y)|<\ep .$$
Hence, $\forall \ep>0,~N>N'$ and $\mu$-a.a. $z\in X$  
$$ \left|H(\frac1N \sum_{k=0}^{N-1} 1_A(f^kz)) - H(\mu(A))\right| < \ep ,$$
which proves that the passage to the limit and the summation can be interchanged.

Now we need to collect together the contributions to the entropy of all the elements of the partition. 
Note that
$$ H(\v{p}(\v{x},n,N)) = \sum_{\Delta_i^n \in \Delta^n} 
                                        H(\frac1{N-n+1} \sum_{k=0}^{N-n} 1_{\Delta_i^n}(f^kz)) .$$
Interchanging the passage to the limit and the finite summation we get
\bea{\a \lim_{n\to\infty}\frac1n\sum_{\Delta_i^n \in \Delta^n} H(\lim_{N\to\infty}\frac1{N-n+1} 
                                      \sum_{k=0}^{N-n} 1_{\Delta_i^n}(f^kz)) \\
\a~~= \lim_{n\to\infty}\lim_{N\to\infty} \frac1n\sum_{\Delta_i^n \in \Delta^n} 
                                   H(\frac1{N-n+1} \sum_{k=0}^{N-n} 1_{\Delta_i^n}(f^kz)) \\
\a~~= \lim_{n\to\infty}\lim_{N\to\infty} \frac1n H(\v{p}(\v{x},n,N)) =  h_{{\rm loc}}(\v{x}|\Delta) ,
}
which yields the result. \qed

\begin{remark}
In particular, in the ergodic case the true limits in the definition of the local entropy exist almost everywhere.
\end{remark}

Theorem~\ref{t:glob-loc-erg} demonstrates that conditional local entropy of a $\mu$-typical 
trajectory of an ergodic dynamical system coincides with the corresponding global conditional 
metric entropy. A similar result cannot be expected in the non-ergodic setting, since 
$\mu$-typical trajectories are extremely rare or even do not exist.

\subsection{Non-ergodic case} \label{s:glob-loc-nonerg}
In this more delicate setting we need to require additionally that $X$ is a Borel subset of a metric compactum.
Let us formulate two results about the ergodic decomposition of invariant measures.

\begin{theorem}[Ergodic Decomposition Theorem, \cite{EW}]\label{t:erg-dec}
For each $f$-invariant measure $\mu$ there exists a set $X^\mu \subset X$ of full $\mu$-measure,  
such that to each $x \in X^\mu$ one can associate a measure $\mu^x$ satisfying 
the following properties:
\begin{enumerate}
\item $\mu^x$ is a $f$-invariant ergodic probability measure on $(X, \Sigma)$.
\item For every $A \in \Sigma$, the map $x \to \mu^x(A)$ is measurable and
   $$
   \mu(A) = \int_X \mu^x(A) \, d\mu(x).
   $$
\item If $\mathcal{E} = \{ E \in \Sigma : f^{-1}E = E \}$ is the $\sigma$-algebra of invariant sets, then there exists a countably generated subalgebra $\tilde{\mathcal{E}} = \sigma(E_1, E_2, \dots) \subset \mathcal{E}$ such that $\mu^x$ is concentrated on the atom
   $$
   [x]_{\tilde{\mathcal{E}}} = \bigcap_{E_i \ni x} E_i \;\cap\; \bigcap_{E_i \not\ni x} E_i^c.
   $$
   Moreover, if $y \in [x]_{\tilde{\mathcal{E}}} \cap X^\mu$, then $\mu^x = \mu^y$.
\item Let a collection of measures $\nu^x$ satisfy the condition 
   $E(\phi|\mathcal{E})(x) = \int \phi(y)~d\nu^x(y)$ for $\mu$-a.a. $x\in X$ and an everywhere dense 
   family of continuous functions $\phi$. Then $\nu^x=\mu^x$ for $\mu$-a.a. $x\in X$.
\end{enumerate}
\end{theorem}

\begin{theorem}[Decomposition of metric entropy, \cite{ELW}] \label{t:ent-dec}
Let $\mu$ be a $f$-invariant measure with ergodic decomposition $\mu(A) = \int_X \mu^x(A) \, d\mu(x)$. Then for any finite measurable partition $\Delta$,
$$
h_\mu(f|\Delta) = \int_X h_{\mu^x}(f|\Delta) \, d\mu(x),
$$
and consequently,
$$
h_\mu(f) = \int_X h_{\mu^x}(f) \, d\mu(x).
$$
\end{theorem}

Now we are ready to formulate and prove our main result about the connection of 
local and global conditional entropies in the non-ergodic case.

\begin{theorem}\label{t:glob-loc-nonerg} 
Let $\mu$ be an $f$-invariant measure. Then 
$$
\int_X h_{\rm{loc}}(\v{x}|\Delta) \, d\mu(x) = h_\mu(f|\Delta) 
$$
for each finite measurable partition $\Delta$. Here we use the notation $x_k:=f^kx~~\forall k\ge0$. 
\end{theorem}

\proof
Take an arbitrary $x \in X^\mu$ (the full-measure set from the ergodic decomposition theorem).  
Applying Theorem~\ref{t:glob-loc-erg} to the space $(X, \Sigma, \mu^x)$ with the ergodic measure 
$\mu^x$ we get that for any finite measurable partition $\Delta$, there exists a set $Z_x$ with $\mu^x(Z_x)=1$ such that for every $y \in Z_x$,
$$
h_{\rm{loc}}(\v{y}|\Delta) = h_{\mu^x}(f|\Delta).
$$

Let us check that 
$h_{\rm{loc}}(\v{x}|\Delta)$ is measurable as a function of $x$.
Indeed,
$$
h_{\rm{loc}}(\v{x}|\Delta) = \limsup_{n \to \infty} ~\limsup_{N \to \infty} \frac{1}{n} 
                                       H\!\left(\bar{p}(\v{x}, f, n, N|\Delta)\right).
$$  
$$
= \limsup_{n \to \infty} ~\limsup_{N \to \infty} \frac{1}{n} \sum_{D \in \Delta^n} H\!\left(\frac{1}{N-n+1} \sum_{j=0}^{N-n} 1_D(f^j x)\right).
$$  
Each indicator function $1_D(f^j x)$ is measurable in $x$; $H$ is continuous, hence measurable; 
finite sums and $\limsup$ preserve measurability. Thus $h_{\rm{loc}}(\cdot|\Delta)$ is measurable. 

Consider a set 
$$
Y = \left\{ y \in X : h_{\mu^y}(f|\Delta) \neq h_{\rm{loc}}(\v{y}|\Delta) \right\}.
$$
By the measurability of the functions $h_{\rm{loc}}(\v{y}|\Delta)$  and  
$y \mapsto h_{\mu^y}(f|\Delta)$ 
(a consequence of the ergodic decomposition theorem), the set $Y$ is measurable.

Let $\t{Y} = Y \cap X^\mu$. For any $x \in X^\mu$, consider $\mu^x(\t{Y})$. 
Because $\mu^x$ is supported on the atom $[x]_{\tilde{\mathcal{E}}}$ (property 3 of the decomposition), and on that atom $\mu^y = \mu^x$ for all $y$ in the atom, we have
$$
\mu^x(\t{Y}) = \mu^x\left( \{ y \in [x]_{\tilde{\mathcal{E}}} \cap X^\mu :
                              ~~ h_{\mu^x}(f|\Delta) \neq h_{\rm{loc}}(\v{y}|\Delta) \} \right).
$$
But by Theorem~\ref{t:glob-loc-erg} applied to the ergodic measure $\mu^x$, the set of $y$ for which $h_{\rm{loc}}(\v{y}|\Delta) \neq h_{\mu^x}(f|\Delta)$ has $\mu^x$-measure zero 
(it is contained in the complement of the set of $\mu^x$-typical points). Hence $\mu^x(\t{Y}) = 0$.

Therefore,
$$
\mu(Y) = \mu(\t{Y} \cap X^\mu) = \mu(\t{Y}) = \int_X \mu^x(\t{Y}) \, d\mu(x) = 0.
$$

Thus, $h_{\rm{loc}}(\v{x}|\Delta) = h_{\mu^x}(f|\Delta)$ for $\mu$-almost every $x$. 
Finally,
$$
h_\mu(f|\Delta) = \int_X h_{\mu^x}(f|\Delta) \, d\mu(x) = \int_X h_{\rm{loc}}(\v{x}|\Delta) \, d\mu(x).
$$
\qed

\section{Connections between metric and local unconditional entropies}\label{s:glob-loc}

\subsection{Family of admissible partitions $\gD^{\v{x}}$}\label{ss:adm-prt}

Both in global and local constructions of dynamical entropies the unconditional versions are 
defined as  supremum of the corresponding conditional versions over certain family of 
partitions. One might expect that in both cases this family can be chosen as all finite 
measurable partitions (as it was done in the Kolmogorov-Sinai construction). 
The following result demonstrates that this is not reasonable for the local setting.

\begin{lemma}\label{l:c-example}
Let $\v{x}$ be an infinite sequence, all elements of which are distinct.
Then $\sup_\Delta h_{{\rm loc}}(\v{x}|\Delta) = \infty$, where the supremum is taken over 
all finite measurable partitions.
\end{lemma}
\proof For each $r\in\IZ_+$ we will construct a measurable partition $\Delta:=\{\Delta_i\}_{i=1}^r$ 
of $(X,\Sigma)$, such that $h_{{\rm loc}}(\v{x}|\Delta)=\log r$. 
By choosing a sequence of values $r$ that increases to infinity, we obtain the result. 
The idea of this construction is to use an $r$-ary sequence $\v{b}:=\{b_i\}_i$ with the maximum 
possible local entropy $\log r$ as the basis for the definition of a special partition $\Delta$, 
which can be done as follows: 
$$ \Delta_j:=\bigcup_{i:b_i=j}\{x_i\}~~j<r, \quad \Delta_r:=X\setminus\bigcup_{j=1}^{r-1}\Delta_j .$$
Obviously, $h_{{\rm loc}}(\v{x}|\Delta) = h_{{\rm loc}}(\v{b})=\log r$. \qed

To demonstrate that a supremum over all finite measurable partitions diverges for an
intuitively ``non-complex'' sequence consider $x_n:=1/n~~n\ge1$. Then $x_n\toas{n\to\infty}0$, 
but $\sup_\Delta h_{\rm loc}(\v{x}|\Delta)=\infty$.

\begin{remark}
Interestingly, the idea used in the proof above does not work for sequences with 
repeated elements, since it does not allow to establish a proper correspondence with 
the chosen $r$-ary sequence $\v{b}$. 
However, there is a more fundamental obstacle.
One might expect that the assumption $x_i\ne x_j~\forall i\ne j$ in Lemma~\ref{l:c-example}
can be easily circumvented by eliminating the extra repeated elements. Indeed, it seems
obvious that the latter operation can only reduce the complexity of the sequence, so the statement
remains true in the case of repeated elements as well.
Unfortunately, this is not true at all. In fact, by inserting new letters into a given sequence,
we may actually decrease its entropy rather than increase it.
Indeed, let $\v{x}$ be a binary sequence with positive local entropy. After inserting
an infinite number of zeros, we obtain a new sequence $\v{y}$ with zero entropy, provided that
the proportion of zeros in $\v{y}$ is large enough (see Lemma~2.4(a) of \cite{Bl25}).
\end{remark}

Lemma~\ref{l:c-example} shows that for a general (say, continuous)
phase space $(X,\Sigma)$ it is necessary to impose certain restrictions on the choice
of the set of admissible partitions (otherwise generic sequences will have infinite entropy).

Observe, that despite the Kolmogorov-Sinai construction do not make explicit 
restrictions of this sort, but implicitly elements of the refined partitions with zero 
measure do not play any role (since $H(0)=0$). Unfortunately this does not help 
much in the case of local entropies and we propose instead a different assumption.

\bdef{We say that a Borel set $A$ is a \emph{set of continuity} of a given measure $\mu$ 
if $\mu(\partial(A))=0$. Here $\partial(A)$ stands for the boundary of $A$.}

\bdef{For a given measure $\mu$ by the set of \emph{admissible partitions}  $\gD^\mu$ we choose 
sets of partitions whose elements belong to the set of continuity for the measure $\mu$.}

One of the main reasons for this choice is explained by the following property, which can be 
found, e.g. in \cite{Bo}.

\begin{lemma}\label{l:weak-conv} 
Let a sequence of measures $\mu_n$ converge weakly to $\mu$. 
Then $\mu_n(A) \toas{n\to\infty}\mu(A)$  for each set $A$ of continuity of $\mu$.
\end{lemma}

On the other hand, one of the main characteristics of a sequence $\v{x}$, which we use is 
a sequence of \emph{empirical} measures $\mu_n^{\v{x}}:=\frac1n\sum_{k=0}^{n-1}\delta_{x_k}$, 
whose weak limits we consider, while the corresponding statistics is evaluated with respect 
to elements of various partitions. 

To be precise from now we assume that the set $X$ is a compact metric space, equipped 
with a metric $\rho$ and a Borel $\sigma$-algebra $\Sigma$.

\bdef{We say that a sequence $\v{x}$ is {\em ergodic} if the measures 
$\mu_n^{\v{x}}\toas{n\to\infty} \mu^{\v{x}}$ weakly.}

Thus we choose $\gD^{\mu^{\v{x}}}$ as the family of admissible partitions for 
an ergodic sequence $\v{x}$. 

For a non-ergodic $\v{x}$ the situation is slightly more complicated. Denote by $\cM(\v{x})$ 
the set of limit points of the sequence $\mu_n^{\v{x}}$. By Prokhorov's Theorem the 
set $\cM(\v{x})$ is compact and non-empty (due to the compactness of $X$). 

For a given limit measure $\nu^{\v{x}} \in \cM(\v{x})$, consider 
$\sup_{\Delta\in \gD^{\nu^{\v{x}}}} h_{{\rm loc}}(\v{x}|\Delta)$. 
By the compactness of $\cM(\v{x})$, the supremum of this quantity over all limit measures 
$\nu^{\v{x}} \in \cM(\v{x})$ is achieved on a certain measure 
$\mu^{\v{x}} \in \cM(\v{x})$ (if it is not unique we choose any of them arbitrarily). 
Finally we give the following definition.

\bdef{For a given sequence $\v{x}\subset X$ by the family of admissible partitions we 
mean $\gD^{\v{x}}:=\gD^{\mu^{\v{x}}}$.}

An important point here is that this family of admissible partitions depends sensitively 
on the sequence $\v{x}$ under consideration.

\subsection{Partitions from $\t\gD^\mu$}\label{ss:partitions}
In the Kolmogorov-Sinai construction the supremum is taken over all measurable partitions,  
therefore one needs to find what happens under the restriction to $\gD^\mu$. 
It turns out that technically it will be simpler for us to deal with a certain subset 
$\t\gD^\mu \subset \gD^\mu$ (which we will construct in this section) and to approximate 
general partitions by this subset later.

Consider all open balls in $X$ of radius $1/n$: $\{B_{1/n}(x)\}_{x\in X}$. 
They form a cover of the compact set $X$. Hence there exists a finite set of points 
$\{x_{n,k}\}_{k=1}^{s(n)}$ such that the balls $\{B_{1/n}(x_{n,k})\}_{k=1}^{s(n)}$ cover $X$. 
Then the set $X^0 := \bigcup_{n=1}^\infty\{x_{n,k}\}_{k=1}^{s(n)}$ is at most countable.

For any $\ep>0$ and any point $x\in X$ there exists $x_{n,k}$ with $\rho(x,x_{n,k})<1/n<\ep$. 
Hence $X^0$ is a dense subset of $X$.

\begin{lemma}\label{l:count-centers}
For a given probability measure $\mu$ and any $x\in X^0$ there exists at most countably many 
$\ep$ for which $\mu(\partial B_\ep(x))>0$.
\end{lemma}
\proof
The sets $\{y\in X \mid \rho(x,y)=\ep\}$ are disjoint for different $\ep$. 
On the other hand, the number of $\ep$ for which $\mu(\{y\in X \mid \rho(x,y)=\ep\}) > 1/m$ is at most $m$, 
because the total measure of $X$ is $1$. 
Hence for each fixed $x\in X^0$ the set of exceptional $\ep$ is at most countable; taking the 
union over the at most countable set $X^0$ we obtain that the set 
$E := \{\ep\in\mathbb{R}_+ \mid \exists x\in X^0:~\mu(\partial B_\ep(x))>0\}$ 
of exceptional $\ep$ is at most countable.

Therefore we can choose $\ep_{x,n}\not\in E$ for each 
$x\in X^0, n\in\mathbb{N}$ such that $\lim_{n\to\infty}\ep_{x,n}=0$, 
$\mu(\partial B_{\ep_{x,n}}(x))=0$ and $\frac1{n+1}<\ep_{x,n}<\frac1n$. \qed

Now we construct a sequence of ``good'' finite partitions $\{\Delta^n\}$ of $X$. 

\begin{lemma}\label{l:good-part}
There exists a sequence $\t\gD^\mu \subset \gD^\mu$ of partitions $\{\Delta^n\}$ such that 
$\bigvee_{n=1}^\infty \Delta^n = \Sigma$.
\end{lemma}
\proof
For each $n$, consider the collection of balls $\{B_{\ep_{x,n}}(x)\}_{x\in X^0}$. 
Since $\ep_{x,n}>1/(n+1)$ and $X^0$ is dense in $X$, these balls cover $X$. Choose a finite subcover 
$\{B_{\ep_{x_i,n}}(x_i)\}_{i=1}^{s(n)}$. From this cover we construct a partition 
$\Delta^n = \{D_n^i\}_{i=1}^{s(n)}$ as follows 
$$
D_n^1 = B_{\ep_{x_1,n}}(x_1),\quad Z_n^1 = D_n^1,\qquad
D_n^i = B_{\ep_{x_i,n}}(x_i) \setminus Z_n^{i-1},\quad Z_n^i = Z_n^{i-1}\cup D_n^i.
$$
Since $\ep_{x,n}<1/n$, the diameter of each $D_n^i$ is less than $2/n$.

Take any open set $A\subset X$. For any $a\in A$ there exists $\delta>0$ such that $B_\delta(a)\subset A$. Choose $N$ with $2/N<\delta$. The diameter of the partition $\bigvee_{n=1}^N \Delta^n$ is less than $2/N<\delta$, so the element of this partition containing $a$ lies inside $B_\delta(a)\subset A$. 
Thus $A$ is a union of elements of the countable collection $\bigcup_{N}\bigvee_{n=1}^N \Delta^n$, 
hence $A\in\bigvee_{n=1}^\infty\Delta^n$. Therefore $\Sigma \subset \bigvee_{n=1}^\infty\Delta^n$. 
The opposite inclusion is obvious, so $\Sigma = \bigvee_{n=1}^\infty\Delta^n$. \qed

\subsection{Refinement of ``good'' partitions}

\begin{lemma}\label{l:L5}
Let $\mu$ be an $f$-invariant probability measure, let $A\in\Sigma$ satisfy $\mu(\partial A)=0$, 
and let $f:(X,\Sigma)\to(X,\Sigma)$ be continuous. 
Then $\mu(\partial(f^{-1}A))=0$.
\end{lemma}
\proof
Continuity of $f$ implies that preimages of open sets are open and preimages of closed sets are closed. 
From ${\rm int}(A)\subseteq A\subset{\rm cl}(A)$ we get
$$
f^{-1}({\rm int}(A)) \subseteq f^{-1}(A) \subseteq f^{-1}({\rm cl}(A)).
$$
Since ${\rm int}(f^{-1}A)$ is the largest open set contained in $f^{-1}A$ and ${\rm cl}(f^{-1}A)$ is the smallest closed set containing $f^{-1}A$,
$$
f^{-1}({\rm int}(A)) \subseteq {\rm int}(f^{-1}A) \subseteq f^{-1}A 
\subseteq {\rm cl}(f^{-1}A) \subseteq f^{-1}({\rm cl}(A)).
$$
Hence
$$
\partial(f^{-1}A) = {\rm cl}(f^{-1}A)\setminus{\rm int}(f^{-1}A) 
     \subseteq f^{-1}({\rm cl}(A)) \setminus f^{-1}({\rm int}(A)) = f^{-1}({\rm cl}(A)\setminus{\rm int}(A)).
$$
Therefore, using the $f$-invariance of $\mu$,
$$
\mu(\partial(f^{-1}A)) \le \mu(f^{-1}(\partial A)) = \mu(\partial A) = 0.
$$
\qed

\begin{lemma}\label{l:L6}
If $A,B\in\Sigma$ satisfy $\mu(\partial A)=\mu(\partial B)=0$, then $\mu(\partial(A\cap B))=0$.
\end{lemma}
\proof 
We have ${\rm int}(A)\subseteq A\subseteq{\rm cl}(A)$ and ${\rm int}(B)\subseteq B\subseteq{\rm cl}(B)$. 
Hence
$$
{\rm int}(A)\cap{\rm int}(B) \subseteq A\cap B \subseteq {\rm cl}(A)\cap{\rm cl}(B).
$$
Moreover,
$$
{\rm int}(A\cap B) \supseteq {\rm int}(A)\cap{\rm int}(B),\qquad
{\rm cl}(A\cap B) \subseteq {\rm cl}(A)\cap{\rm cl}(B).
$$
Thus
$$
\partial(A\cap B) = {\rm cl}(A\cap B)\setminus{\rm int}(A\cap B) 
                      \subseteq \bigl({\rm cl}(A)\cap{\rm cl}(B)\bigr) 
                            \setminus \bigl({\rm int}(A)\cap{\rm int}(B)\bigr) 
                      \subseteq \partial A \cup \partial B.
$$
Consequently,
$$
\mu(\partial(A\cap B)) \le \mu(\partial A) + \mu(\partial B) = 0,
$$
which finishes the proof. \qed

\begin{lemma}\label{l:subpart-cont}
All subpartitions of the partitions constructed in the proof of Lemma~\ref{l:good-part} are sets of 
continuity of the measure $\mu$.
\end{lemma}
\proof
A set and its complement have the same boundary, so if a set has boundary of $\mu$-zero measure 
then its complement also has $\mu$-zero boundary. The difference $A\setminus B = A\cap B^c$; 
if $A$ and $B$ have $\mu$-zero boundary, then by Lemma~\ref{l:L6} (applied twice) $A\setminus B$ 
also has $\mu$-zero boundary. Similarly, $A\cup B$ has $\mu$-zero boundary because 
$(A\cup B)^c = A^c\cap B^c$. The partitions $\Delta^n$ are constructed from balls 
$B_{\ep_{x_i,n}}(x_i)$ which have zero boundary by construction, using only finite unions, intersections and differences. Hence each $\Delta^n$ consists of sets with $\mu$-zero boundary. Any subpartition $\Delta\subset\Delta^n$ is a finite union of elements of $\Delta^n$, therefore also 
has $\mu$-zero boundary.
\qed

\subsection{Ergodic case}

The following statement collects a number of well known results about ergodic invariant measures. 

\begin{lemma}\label{l:erg-mes} (See, e.g. \cite{Bl17, Ba}) 
Denote $Z_\mu:=\{x\in X:~ \frac1n\sum_{k=0}^{n-1} \delta_{f^kx} \toas{n\to\infty}\mu\}$, 
${\rm Erg}(f)$ - the set of all ergodic $f$-invariant measures, $Z(f):=\cup_{\mu\in {\rm Erg}(f)}Z_\mu$. 
Then $Z_\mu, Z(f)$ are Borel sets, and $\mu(Z_\mu) = \mu(Z(f))=1~~\forall \mu\in {\rm Erg}(f)$.
\end{lemma}

Now we are ready to study the connection between local and metric entropy.

\begin{theorem}\label{t:global-erg}
Let $\mu$ be a probability ergodic measure for a continuous transformation $f:(X,\Sigma)\to(X,\Sigma)$.  Then there exists a set $S\subset X$ of full $\mu$-measure such that for every $x\in S$,
$$
h_{{\rm loc}}(\v{x}) = h_\mu(f).
$$
\end{theorem}
\proof
Set $\tilde\Delta^n := \bigvee_{i=1}^n \Delta^i$. Then 
$\tilde\Delta^1\subset\tilde\Delta^2\subset\cdots$. Moreover,
$$
\Delta^n \subset \tilde\Delta^n \subset \bigvee_{i=0}^\infty f^{-i}\tilde\Delta^n,
$$
and
$$
\Sigma = \bigvee_{n=1}^\infty \Delta^n 
       \subset \bigvee_{n=1}^\infty \bigvee_{i=0}^\infty f^{-i}\tilde\Delta^n \subset \Sigma,
$$
so the property \ref{p:entropy}(2) applies:
$$
h_\mu(f) = \lim_{n\to+\infty}\sup_{\Delta\subset\tilde\Delta^n} h_\mu(f|\Delta).
$$

Consequently,
$$
h_\mu(f) = \lim_{n\to+\infty}\sup_{\Delta\subset\tilde\Delta^n} h_\mu(f|\Delta) 
              = \sup_{\Delta\in\t\gD^\mu} h_\mu(f|\Delta).
$$
From the weak convergence of ergodic averages to $\mu$, it follows that Birkhoff's theorem
is true on $Z_\mu$ simultaneously for all continuity sets, and therefore
from the proof of Theorem~\ref{t:glob-loc-erg} we deduce that 
$h_{\rm loc}(\v{x}|\Delta) = h_\mu(f|\Delta)$ for each $x\in Z_\mu$ and 
any partition $\Delta\in \t\gD^\mu$. 

Taking supremum over $\Delta\in\t\gD^\mu$ gives
$$
h_\mu(f) = h_{{\rm loc}}(\v{x}).
$$
\qed

\subsection{Non-ergodic case}

\begin{theorem}\label{t:global-nonerg}
Let $\mu$ be a probability invariant measure for a continuous transformation $f:(X,\Sigma)\to(X,\Sigma)$. Then 
$$
\int_X h_{{\rm loc}}(\v{x})\,d\mu(x) = h_\mu(f).
$$
\end{theorem}

Before giving a proof we need a few technical results.

First, let us demonstrate that for a measure preserving dynamical system 
$(f,X,\mu)$ the set of untypical points is negligible.

\begin{lemma} \label{l:emp-to-erg}
The set of points $x$ for which the empirical measures 
$m_n^x:=\frac1n \sum_{k=0}^{n-1}\delta_{f^kx}$ 
converge weakly  to a non-ergodic $f$-invariant measure
is a $\mu$-null set for any $f$-invariant measure $\mu$.
\end{lemma}

\proof Follows from Lemma~\ref{l:erg-mes}. \qed

Fix an $f$-invariant measure $\mu$. By Lemmas~\ref{l:erg-mes} and \ref{l:emp-to-erg} for $\mu$-a.a. 
points $x\in X$ the sequence of empirical measures $m_n^x$ converges weakly to some 
ergodic measure $\nu^x$.

\begin{lemma} \label{l:l11}
$\nu^x = \mu^x$ for $\mu$-a.a. points $x\in X$, where $\mu^x$ are 
measures from the ergodic decomposition of $\mu$.
\end{lemma}
\proof By the definition of weak convergence for each continuous function $\phi$ we have 
$$ \lim_{n\to\infty} \frac1n \sum_{k=0}^{n-1}\phi(f^kx) = \int \phi(y) ~d\nu^x(y). $$
On the other hand, rewriting Birkhoff's Ergodic Theorem in terms of conditional 
mathematical expectations, we get
$$ \lim_{n\to\infty} \frac1n \sum_{k=0}^{n-1}\phi(f^kx) = E(\phi|\mathcal{E})(x) $$
$\mu$-a.e. for a dense family of continuous functions $\phi$. The result follows from 
the 4-th item of Theorem~\ref{t:erg-dec}. \qed

\bigskip

\n
{\bf Proof} of Theorem~\ref{t:global-nonerg}. 
Applying Theorem~\ref{t:global-erg} to the ergodic dynamical system $(f,X,\nu^x)$ 
we get 
$$ h_{\rm loc}(\v{y}) = h_{\nu^x}(f) = h_{\nu^y}(f) $$
for each point $y\in Z_{\nu^x}$. 

Therefore $h_{\rm loc}(\v{y}) = h_{\nu^y}(f)$ for each $y\in \cup_{\nu\in{\rm Erg(f)}}Z_\nu$. 
Since $\mu(Z(f))=1$ (see Lemma~\ref{l:erg-mes} together with the ergodic decomposition), 
we have $h_{{\rm loc}}(\v{y}) = h_{\nu^y}(f)$ 
for $\mu$-a.a. $y\in X$. 

Now by Lemma~\ref{l:l11} we know that $\nu^x = \mu^x$ for $\mu$-a.a. points $x\in X$, 
hence $h_{\nu^x}(f) = h_{\mu^x}(f)$ for $\mu$-a.a. $x$, which implies 
that $h_{{\rm loc}}(\v{x}) = h_{\mu^x}(f)$ for $\mu$-a.a. points $x\in X$.

Finally the entropy decomposition Theorem~\ref{t:ent-dec} gives 
$$
h_\mu(f) = \int_X h_{\mu^x}(f)\,d\mu(x) = \int_X h_{{\rm loc}}(\v{x})\,d\mu(x).
$$
\qed

\section{Connections between topological and information entropies}\label{s:top-Inf}

First of all, let us explain why one cannot use coverings instead of partitions in the 
construction of the information entropy. The main problem here is that the elements 
of the covering intersect between themselves, and hence some points belong to several 
elements of the partition and thus get multiple codings. 

Consider an example. Let $x_i\equiv x~~\forall i$. Then the 
corresponding sequence cannot be considered to be complex. On the other 
hand we will show that for each $r\in\IZ_+$ one can construct a collection 
of $r$ non-coinciding open sets $\Delta^{(r)}$ such that the point $x$ belongs 
to each element of $\Delta^{(r)}$. Using this coding, we get $L(\v{x},n,N)=r^n$. 
Therefore the corresponding conditional information entropy is equal to $r$. 
Since $r$ is arbitrary, this leads to the infinite information entropy.

\bigskip

Thus, being unable to calculate the information entropy via coverings, 
in order to compare the topological entropy of a dynamical system with 
information entropies of its individual trajectories one needs to formulate 
a version of the topological entropy via partitions. 

Let the map $f:X\to X$ defines a dynamical system and let $\Delta$ be a finite 
partition of $X$. 
Introduce a transition matrix $M:=\{m_{ij}\}$, where $m_{ij}=1$ 
if $\Delta_i \cap f^{-1}\Delta_j \ne \emptyset$ and $=0$ otherwise.

Then on the Cantor set $X_M$ (the space of sequences with 
the alphabet $\cA:=\{1,2,\dots,r\}$ with the  transition matrix $M$)  
the left shift map $\sigma$ defines a symbolic dynamical system with trajectories of type: 

$$ \v{x} = (x_1x_2\dots x_k\underbrace{x_{k+1}\dots x_{k+n}}_{\v{w}} 
                 x_{k+n+1}\dots  x_N\dots) $$
Denoting by $A_\Delta^n$ the set of all admissible words $\v{w}$ of length $n$  
(i.e. different pieces of length $n$ of all trajectories of $(\sigma,X_M)$)  
and by $\#A_\Delta^n$ -- the number of such words, we set  
$$h_{{\rm top}}^*(f|\Delta):=\lim_{n\to\infty}\frac1n \log(\#A_\Delta^n) .$$
Sometimes this quantity is referred as a conditional {\em combinatorial} entropy. 

\begin{lemma}\label{l:top-info}
For each trajectory $\v{x}$ of the map $f$ and any finite partition $\Delta$ 
$$ h_{{\rm top}}^*(f|\Delta) \ge h_{{\rm info}}(\v{x}|\Delta) .$$
\end{lemma}
\proof Indeed, all admissible words appearing in $\v{x}$ are taken into account in 
the definition of the combinatorial entropy. \qed

It remains to find out the connection between the combinatorial and topological entropies. 
In general the answer is not known, but according to \cite{KH} at least for 
expansive systems $h_{top}(U) \le h_\mu(f|\Delta)$ for all open covers $U$,  
some ergodic $\mu$ and each partition $\Delta$ finer than $U$.

\bigskip

In the case of unconditional information entropy, the same problem of choosing admissible 
partitions arises as in the case of local entropy. Indeed, by Lemma~\ref{l:ineq-entr} 
$$ h_{{\rm info}}(\v{x}|\Delta) \ge h_{{\rm loc}}(\v{x}|\Delta)  $$
for any pair of $\v{x}, \Delta$. Thus the claim of Lemma~\ref{l:c-example} remains 
valid for the information entropy as well. 

Therefore to define the unconditional information entropy we use the same family of 
admissible partitions introduced in Section~\ref{ss:adm-prt}:
$$ h_{{\rm info}}(\v{x}) := \sup_{\Delta\in \gD^{\v{x}}}h_{{\rm info}}(\v{x}|\Delta) .$$

This allows to extend the result of Lemma~\ref{l:ineq-entr} for unconditional entropies.

\begin{lemma}\label{l:loc-info-uncond} 
$$ h_{{\rm loc}}(\v{x}) \le h_{{\rm info}}(\v{x}) .$$   
\end{lemma}

\section{Topological invariance of local dynamical entropies}\label{s:topology}

\begin{theorem}\label{t:top-inv}
Let $\v\alpha = \alpha_0, \alpha_1, \dots$ and $\v\beta = \beta_0, \beta_1, \dots$ be sequences of points 
in metric compacts $A$ and $B$, respectively. If there exist two compact subsets $A' \subset A$ and 
$B' \subset B$ containing all elements of the sequences $\alpha$ and $\beta$, respectively, and a homeomorphism $\varphi : A' \to B'$ mapping the sequence $\v\alpha$ element-by-element onto 
the sequence $\v\beta$, then the unconditional local entropies of the sequences $\v\alpha$ and $\v\beta$ coincide.
\end{theorem}

\begin{proof}
First, note that the local entropy of the sequence $\v\alpha$ does not change when we pass from the 
compact $A$ to its compact subset $A'$ (and similarly for the other sequence). Indeed, consider an arbitrary 
limit measure $\mu$ of the sequence $\frac{1}{n}\sum_{i=0}^{n-1} \delta_{\alpha_i}$. All measures $\frac{1}{n}\sum_{i=0}^{n-1} \delta_{\alpha_i}$ are measures on $A'$, since the sequence $\v\alpha$ 
lies in $A'$. Consequently, the partial limits of the sequence 
$\frac{1}{n}\sum_{i=0}^{n-1} \delta_{\alpha_i}$ 
lie in the compact space of probability measures on $A'$ (which is compact since $A'$ is compact). 
Thus all limit measures of the sequence $\frac{1}{n}\sum_{i=0}^{n-1} \delta_{\alpha_i}$ in the 
compact space of probability measures on $A$ are in fact probability measures on $A'$. 
In other words, for every limit measure $\mu$, we have $\mu(A \setminus A') = 0$.

Thus the set of measures used in the definition of the entropy of the sequence $\v\alpha$ on $A$ and 
on $A'$ coincide. Fix an arbitrary measure $\mu$ from this set. 


Let $\Delta = \{D_i\}_{i=1}^r$ be a partition of $A$. Associate to it the collection 
$\Delta' = \{D'_i\}_{i=1}^r$ as follows: $D'_i = A' \cap D_i$. This collection obviously consists of disjoint 
sets and covers all of $A'$, since $A' \subset A$; hence it is a partition. Let us check that the 
$\mu$-measure of boundary of the sets in this collection (as subsets of $A'$, not $A$) is zero. 
To do this, it suffices to observe that 
$\partial_{A'}(D_i\cap A') \subseteq \partial_A(D_i)\cap A'$; indeed, a point 
$a\in A'$ that is interior for $D_i$ in $A$ (i.e.\ some $A$-neighbourhood of $a$ lies in $D_i$) 
is interior for $D_i\cap A'$ in $A'$, and similarly for the complement. 
Consequently $\mu(\partial_{A'}(D'_i)) \le \mu(\partial_A(D_i)) = 0$, and no new boundary points appear.

Now suppose we are given a partition $\Delta' = \{D'_i\}_{i=1}^r$ of $A'$ with 
$\mu$-zero boundary. Construct from it a collection $\Delta = \{D_i\}_{i=1}^r$ 
as follows: each set $D'_i$ consists of its interior points 
${\rm int}(D'_i)$ and boundary points $D'_i \setminus {\rm int}(D'_i)$ (considering $D'_i$ as a subset of $A'$); the first is measurable since it is open, and the second as the difference of two measurable sets. To each interior point $x \in {\rm int}(D'_i)$ assign the open ball in $A$ of radius half the distance to the closure of the complement of $D'_i$ in $A'$, i.e. 
$B_{\frac{1}{2} \rho(x, (A' \setminus D'_i))}(x)$. This definition is correct since $(A' \setminus D'_i)$ 
is a closed subset of a compact, hence compact.

The desired set is
$$ D_i = (D'_i \setminus {\rm int}(D'_i)) \cup 
     \bigcup_{x \in {\rm int}(D'_i)} B_{\frac{1}{2} \rho(x, (A' \setminus D'_i))}(x). $$
The sets $D_i$ are measurable, since they are the union of the measurable set $D'_i \setminus {\rm int}(D'_i)$ and the open set 
$\bigcup_{x \in {\rm int}(D'_i)} B_{\frac{1}{2} \rho(x, (A' \setminus D'_i))}(x)$.

The sets $D_i$ are disjoint. Indeed, suppose $x \in D_i$ and $x \in D_j$ with $i \neq j$. Then there are three different cases.

\begin{itemize}
\item $x \in D'_i \setminus {\rm int}(D'_i)$ and $x \in D'_j \setminus {\rm int}(D'_j)$, but then $x \in D'_i \cap D'_j$, which is impossible since $D'_i \cap D'_j = \varnothing$.

\item $x \in D'_i \setminus {\rm int}(D'_i)$ and 
$x \in B_{\frac{1}{2} \rho(y, (A' \setminus D'_j))}(y)$. Note that $x$ is contained in $D'_i \subset A' \setminus D'_j$, hence $x \in (A' \setminus D'_j)$. Thus 
$\rho(x, y) \ge \rho(y, (A' \setminus D'_j)) \ge \frac{1}{2} \rho(y, (A' \setminus D'_j))$. 
Consequently, $x \notin B_{\frac{1}{2} \rho(y, (A' \setminus D'_j))}(y)$. Contradiction.

\item $x \in B_{\frac{1}{2} \rho(y, (A' \setminus D'_i))}(y) ~\cap~ 
                 B_{\frac{1}{2} \rho(z, (A' \setminus D'_j))}(z)$. 
Then
$$
\rho(y, z) \le \rho(y, x) + \rho(x, z) < \frac{1}{2}\bigl(\rho(y, (A' \setminus D'_i)) + \rho(z, (A' \setminus D'_j))\bigr).
$$
On the other hand, $y \in {\rm int}(D'_i) \subset (A' \setminus D'_j)$ and $z \in {\rm int}(D'_j) \subset (A' \setminus D'_i)$, hence
$$
\frac{1}{2}\bigl(\rho(y, (A' \setminus D'_i)) + \rho(z, (A' \setminus D'_j))\bigr) \le \frac{1}{2}(\rho(y, z) + \rho(z, y)) = \rho(y, z).
$$
Contradiction.
\end{itemize}

Thus, to obtain a measurable partition of $A$ from the collection $\{D_i\}_{i=1}^r$, we must add 
the set $D_0 = A \setminus (\bigcup_{i=1}^r D_i)$.

Let us check that the sets of this partition have $\mu$-zero boundary. 
For $i \ge 1$, consider an arbitrary boundary point $x$ of the set $D_i$ such that $x \in A'$. 
Note that $D_i \cap A' = D'_i$, since by definition $D'_i \subset D_i$ and $D_i \cap D_j = \emptyset$ 
for $i \neq j$. Thus either $x \in \partial(D'_i)$ or $x \in {\rm int}(D'_i)$; but if 
$x \in {\rm int}(D'_i)$, then $x$ is contained in $D_i$ together with a ball of positive radius 
$\frac{1}{2} \rho(x, A' \setminus D'_i)$. Thus every boundary point of $D_i$ lies either outside $A'$ 
or on the boundary of $D'_i$. Consequently, the measure of the boundary remains zero, since 
$\mu(A \setminus A') = 0$.

The boundary of the set $D_0$ coincides with the boundary of the set $\bigcup_{i=1}^r D_i$, which 
is of zero measure by the previously proved assertion, since $\mu(\partial(D_i)) = 0$.

Thus, for every partition of $A$ we have constructed a partition of $A'$ indistinguishable from the point 
of view of the sequence $\v\alpha$, and conversely. Consequently, we may assume that our original 
spaces $A$ and $B$ coincide with $A'$ and $B'$ from the condition of the theorem.

Note that the image of the delta measure at a point $a \in A$ under the transformation $\varphi$ is 
the delta measure at the point $\varphi(a) \in B$. Indeed, $\varphi(a) \in C$ if and only if 
$a \in \varphi^{-1}C$; in other words, $\delta_{\varphi(a)}(C) = \delta_a(\varphi^{-1}C)$.

Consider the sequence of empirical measures 
$\delta_n^{\v\alpha} := \frac{1}{n}\sum_{i=0}^{n-1} \delta_{\alpha_i}$, 
corresponding to the sequence $\v\alpha$. Let $\lim_{k\to\infty} \delta_{n_k}^{\v\alpha} = \mu$, 
and let $\varphi\mu$ be the image of the measure $\mu$ under the action of the map $\varphi$, 
while $\varphi(\v\alpha)$ is the sequence $\varphi(\alpha_i)$.
Let us check that then $\lim_{k\to\infty} \delta_{n_k}^{\varphi(\v\alpha)} = \varphi\mu$. 
Fix an arbitrary continuous function $g : B \to \mathbb{R}$ and consider the function 
$g(\varphi(\cdot)) : A \to \mathbb{R}$. It is continuous, since $g$ and $\varphi$ are continuous, 
and bounded, as a continuous function on a compact. Then by the definition of weak convergence for the continuous bounded function $g(\varphi(\cdot))$,
$$
\lim_{k \to +\infty} \int_A g(\varphi)\, d\delta_{n_k}^{\v\alpha} = \int_A g(\varphi)\, d\mu.
$$
Passing in the integrals from the measure to its image,
$$
\lim_{k \to +\infty} \int_B g\, d\delta_{n_k}^{\varphi(\v\alpha)} = \int_B g\, d\varphi\mu.
$$
Since $\varphi(\v\alpha) = \beta$, the sequence $\delta_{n_k}^{\v\beta}$ weakly converges to $\varphi\mu$.

Consequently, we have obtained that if the sequence $\delta_n^{\v\alpha}$ has a convergent 
subsequence, then the corresponding subsequence of the sequence $\delta_n^{\v\beta}$ also converges. Since $\varphi$ is a homeomorphism, the converse is also true. In other words, the limit points of the sequences of measures $\delta_n^{\v\alpha}$ and $\delta_n^{\v\beta}$ are in one-to-one 
correspondence given by the transformation $\varphi$.

Consider an arbitrary partition $\Delta := \{D_i\}_{i=1}^r$ of $A$. It corresponds one-to-one to 
the partition $\varphi(\Delta) := \{\varphi(D_i)\}_{i=1}^r$. The condition $\alpha_i \in D_j$ 
is equivalent to the condition $\beta_i = \varphi(\alpha_i) \in \varphi(D_j)$; consequently, 
the sequences of indices assigned in the definition of conditional local entropy to the sequences 
$\v\alpha$ and $\v\beta$ will be the same. Thus
$$
h_{\rm loc}(\v\alpha|\Delta) = h_{\rm loc}(\v\beta|\varphi(\Delta)).
$$
Note that since $\varphi$ is a homeomorphism and a homeomorphism preserves boundaries, for 
every limit measure $\mu$ we have $\mu(\partial D_i) = \varphi\mu(\partial \varphi(D_i))$ for every $i$. 
That is, if the partition $\Delta$ had zero boundary with respect to $\mu$, then its image will 
have zero boundary with respect to $\varphi\mu$. In other words, partitions of $A$ with zero boundary 
with respect to $\mu$ are in one-to-one correspondence with partitions of $B$ with zero boundary 
with respect to $\varphi\mu$; moreover, this correspondence preserves conditional local entropy.

Thus, if for a fixed limit measure $\mu$ we consider the supremum of the conditional local entropy 
of the sequence $\v\alpha$ over all partitions with $\mu$-zero boundary, we obtain the same result 
as taking the supremum of the conditional local entropy of the sequence $\v\beta$ over all partitions with $\varphi\mu$-zero boundary.

Finally, since $\varphi$ also gives a one-to-one correspondence of the limit measures of the sequences $\v\alpha$ and $\v\beta$, taking the supremum of the unconditional local entropies over all limit 
measures for the sequences $\v\alpha$ and $\v\beta$, we again obtain the same results. \qed
\end{proof} 



\newpage
\begin{thebibliography}{99} \label{bibl}

\bibitem{AKM} R.L. Adler, A.G. Konheim, M.H. McAndrew. 
   {\em Topological Entropy}. Transactions of the American Mathematical Society. 114:2 (1965): 309-319. doi:10.1090/S0002-9947-1965-0175106-9

\bibitem{Ba} V.I. Bakhtin, 
  {\em Ergodic Support of a Dynamical System and a Natural Representation of Choquet Distributions 
  for Invariant Measures}, Math Notes 120 (2026), 392-395. https://doi.org/10.1134/S0001434626603436 

\bibitem{Bi} P. Bilingsley. {\em Ergodic Theory and Information}, Wiley \& Sons, 1965.


\bibitem{Bl17} M. Blank. {\em Ergodic averaging with and without invariant measures}, 
   Nonlinearity, 30:12(2017), 4649-4664. DOI: 10.1088/1361-6544/aa8fe8     

\bibitem{Bl25} M. Blank. {\em Are prime numbers and quadratic residues random?}
  Discrete and Continuous Dynamical Systems, 45:11(2025), 4186-4201.  Doi: 10.3934/dcds.2025053 


\bibitem{Bo} V.I. Bogachev. {\em Weak convergence of measures}, 
 AMS, Mathematical Surveys and Monographs, V.234; 2018; 286 pp. ISBN: 978-1-4704-7798-1

\bibitem{BK} M. Brin and A. Katok, {\em On local entropy}, in Geometric dynamics 
   (Rio de Janeiro, 1981), 30-38, Springer, Berlin, 1983.

\bibitem{Br} A. A. Brudno. 
  {\em The complexity of the trajectories of a dynamical system}, 
  Russian Math. Surveys, 33:1 (1978), 197-198.

\bibitem{EW} M. Einsiedler, T. Ward. 
  {\em Ergodic Theory with a view towards Number Theory}, 
  Springer London, 2011 DOI 10.1007/978-0-85729-021-2

\bibitem{ELW} M. Einsiedler, E. Lindenstrauss, Th. Ward, \emph{Entropy in Ergodic Theory and Topological Dynamics}, Chapter 2, pp. 36-37. (In preparation.) 
https://tbward0.wixsite.com/books/entropy

\bibitem{Ka} A. Katok. 
  {\em Fifty years of entropy in dynamics: 1958-2007}, 
  Journal of Modern Dynamics,  1: 4 (2007), 545-596.  Doi: 10.3934/jmd.2007.1.545

\bibitem{KH}  A. Katok, B. Hasselblatt, Introduction to the Modern Theory of Dynamical Systems,
Cambridge University Press, 1995.

\bibitem{Wa} P. Walters. 
  {\em An Introduction to Ergodic Theory}, Springer New York, NY, 1982. 

\end{thebibliography}
\end{document}